\documentclass[12pt,a4paper,reqno]{article}%
\usepackage{amsfonts}
\usepackage{amsmath}
\usepackage{amssymb}
\usepackage[all]{xy}
\usepackage{rotating}
\usepackage{graphicx}
\DeclareFontFamily{OT1}{wncyi}{}
\DeclareFontShape{OT1}{wncyi}{m}{it}{
<5> <6> <7> <8> <9> gen * wncyi
<10> <10.95> <12> <14.4> <17.28> <20.74> <24.88> wncyi10
}{}
\DeclareSymbolFont{cyrletters}{OT1}{wncyi}{m}{it}
\DeclareSymbolFontAlphabet{\cyrmath}{cyrletters}
\DeclareMathSymbol{\rE}{\cyrmath}{cyrletters}{003}
\newtheorem{theorem}{Theorem}
\newtheorem{corollary}[theorem]{Corollary}
\newtheorem{definition}[theorem]{Definition}
\newtheorem{example}[theorem]{Example}

\newtheorem{proposition}[theorem]{Proposition}
\newtheorem{remark}[theorem]{Remark}

\title{\bf Iterated Differential Forms V:  $\mathcal{C}$--Spectral Sequence on $J^{\infty
}(\pi)$}
\author{\sc{A.~M.~Vinogradov}\thanks{{\bf e}-{\it mail}: \texttt{vinograd@unisa.it}} and \sc{L.~Vitagliano}\thanks{{\bf e}-{\it mail}: \texttt{luca\_vitagliano@fastwebnet.it}}\\
\small{DMI, Universit\`a degli Studi di Salerno}\\ \small{and INFN, Gruppo collegato di Salerno,}\\
\small{Via Ponte don Melillo, 84084 Fisciano (SA), Italy}}
\begin{document}
\maketitle
\begin{abstract}
In the preceding note \cite{vv07} the
$\Lambda_{k-1}\mathcal{C}$--spectral sequence, whose first term is
composed of \emph{secondary iterated differential forms}, was
constructed for a generic diffiety. In this note the zero and
first terms of this spectral sequence are explicitly computed for
infinite jet spaces. In particular, this gives an explicit
description of secondary covariant tensors on these spaces and
some basic operations with them. On the basis of these results a
description of the $\Lambda_{k-1}\mathcal{C}$--spectral sequence
for infinitely prolonged PDE's will be given in the subsequent
note.
\end{abstract}

\maketitle
\newpage
Introduced in \cite{vv07} secondary iterated
differential forms on a generic diffiety
$(\mathcal{O},\mathcal{C})$ are elements of the first term of the
$\Lambda_{k-1}\mathcal{C}$-- spectral sequence associated with it.
In the present note we give an explicit description of the zeroth
and first terms of this spectral sequence for the infinite jet
space of a vector bundle $\mathcal{O}=J^{\infty}(\pi)$,
$\pi:E\longrightarrow M$. Moreover, analogues of basic operations
of tensor analysis in secondary calculus on these spaces are
constructed. This goal is got by a due generalization of the
approach developed in \cite{v84} taking into account some
improvements proposed in \cite{kv98}.
\section{Notations and conventions}

In what concerns geometry of infinite jet manifolds and general
theory of iterated differential forms (shortly, IDFs), we follow
\cite{b99} and \cite{vv06}, respectively. The notation is slightly
simplified against \cite{b99,vv06}. Namely, when the context
allows we omit the reference to arguments \textquotedblleft$\pi$\textquotedblright\ and
\textquotedblleft$J^{\infty}(\pi)$\textquotedblright. For instance, we use $\mathcal{F}$ and
$\Lambda$ for algebras of smooth functions and differential forms
on $J^{\infty}(\pi)$, respectively, instead of standard
$\mathcal{F}(\pi)$ and $\Lambda(J^{\infty}(\pi))$, etc.

Also the algebra of horizontal differential forms on
$J^{\infty}(\pi)$ is denoted by $\mathcal{H} \Lambda$ (rather than
$\Lambda_{0}$) and the algebra of vertical differential forms by
$\mathcal{C}^{\bullet}\Lambda$. Accordingly, $d^{\,\mathsf{h}}$
and $d^{\,\mathsf{v}}$ stand for horizontal and vertical
differentials, respectively.

As usually, $\{\ldots,x^{\mu},\ldots,u^{j},\ldots\}$, $\mu=1,\ldots,n$, $j=1,\ldots,m$, $n=\dim M$, $m=\dim E-n$, denote an adapted to $\pi$ local
chart on $E$, assuming that $\{\ldots,x^{\mu},\ldots\}$ is a local chart on $M$.
We follow the notation of \cite{vv07} in what concerns the theory
of secondary IDF's. The wedge product of IDFs is denoted,
according to \cite{vv07}, by the \textquotedblleft dot\textquotedblright. Finally, the fact that an
isomorphism is \emph{canonical} is stressed using the symbol
$\simeq$.

\section{IDF\lowercase{s} on $J^{\infty}(\pi)$}

Recall that, due to the fibered structure of $J^{\infty}(\pi)$,
the de Rham differential $d$ splits into vertical and horizontal
parts: $d=d^{\,\mathsf{v}}+d^{\,\mathsf{h}}$. $(\Lambda,d^{\,\mathsf{v}%
},d^{\,\mathsf{h}})$ is a bi--complex, for the first time
introduced in \cite{v78} (see also \cite{t80,t82,v84}) and later called the
\emph{variational bicomplex} (see \cite{kv98,b99}, see also \cite{a92}). A variational
bi--complex exists locally for any
diffiety and may be understood as a local description of the $\mathcal{C}%
$--spectral sequence. Moreover, $\Lambda^{1}$ splits naturally as
$\Lambda ^{1}=\mathcal{C}\Lambda^{1}\oplus\mathcal{H}\Lambda^{1}$.
Accordingly, $\Lambda$ is factorized as
$\Lambda\simeq\mathcal{C}^{\bullet}\Lambda
\otimes_{\mathcal{F}}\mathcal{H}\Lambda$
($\mathcal{C}^{\bullet}\Lambda
\equiv\bigwedge^{\bullet}\mathcal{C}\Lambda^{1}$,
$\mathcal{H}\Lambda
\equiv\bigwedge^{\bullet}\mathcal{H}\Lambda^{1}$). In a dual
manner,
$\mathrm{D}$ splits as $\mathrm{D}=\mathcal{C}\mathrm{D}\oplus\mathrm{D}%
^{\,\mathsf{v}}$ (see \cite{b99}).

Note that $d^{\,\mathsf{v}}$ and $d^{\,\mathsf{h}}$ extends as
derivations to $\Lambda_{k}$. Abusing the notation we continue to
denote  these extensions by $d^{\,\mathsf{v}}$ and
$d^{\,\mathsf{h}}$. Put $d_{m}^{\,\mathsf{v}}=\kappa _{1m}\circ
d^{\,\mathsf{v}}\circ\kappa_{1m}:\Lambda_{k}\longrightarrow
\Lambda_{k}$, \; $d_{m}^{\,\mathsf{h}}=\kappa_{1m}\circ d^{\,\mathsf{h}}%
\circ\kappa_{1m}:\Lambda_{k}\longrightarrow\Lambda_{k}$, $\kappa_{1m}%
:\Lambda_{k}\longrightarrow\Lambda_{k}$ being the involution that interchanges
$d_{1}$ and $d_{m}$, $m\leq k$ and also $d_{K}^{\,\mathsf{v}}=d_{k_{1}%
}^{\,\mathsf{v}}\circ\cdots\circ d_{k_{s}}^{\,\mathsf{v}}$, $K=\{k_{1}%
,\ldots,k_{s}\}\subset\{1,\ldots,k\}$. $(\Lambda_{k},(d_{1}^{\,\mathsf{v}%
},d_{1}^{\,\mathsf{h}}),\ldots,(d_{k}^{\,\mathsf{v}},d_{k}^{\,\mathsf{h}}))$,
where $d_{1}^{\,\mathsf{v}}=d^{\,\mathsf{v}}$ and $d_{1}^{\,\mathsf{h}%
}=d^{\,\mathsf{h}}$, is a multiple bi--complex. Define $\mathcal{H}\Lambda
_{k}^{1}\subset\Lambda_{k}^{1}$ as the $\Lambda_{k-1}$--submodule generated by
$\pi_{\infty}^{\ast}(\Lambda_{k}^{1}(M))$, $\pi_{\infty}:J^{\infty}%
(\pi)\longrightarrow M$ being a natural projection. $\mathcal{C}\Lambda
_{k}^{1}\subset\Lambda_{k}^{1}$ (see \cite{vv07}) is a locally free
$\Lambda_{k-1}$--submodule generated by $d_{k}^{\,\mathsf{v}}(\Lambda_{k-1})$.
Elements
\begin{equation}
\mathcal{C}\Lambda_{k}^{1}\ni d_{k}^{\,\mathsf{v}}d_{K}^{\,\mathsf{v}%
}u_{\sigma}^{j},\quad K\subset\{1,\ldots,k-1\},\;j=1,\ldots,m,\;\sigma\text{ a
lenght $n$ multi--index} \label{basis1}%
\end{equation}
form a local basis of it and must be understood as iterated
Cartan forms. Similarly,
$\mathcal{H}\Lambda_{k}^{1}\subset\Lambda_{k}^{1}$ is a locally
free $\Lambda_{k-1}$--submodule generated by $d_{k}^{\,\mathsf{h}}%
(\Lambda_{k-1})$ and elements
\begin{equation}
\mathcal{H}\Lambda_{k}^{1}\ni d_{k}d_{K}x^{\mu},\quad K\subset\{1,\ldots
,k-1\},\;\mu=1,\ldots,n \label{basis2}%
\end{equation}
form an its local basis.

Denote by $\mathcal{C}^{\bullet}\Lambda_{k}\subset\Lambda_{k}$
(resp.~$\mathcal{H}\Lambda_{k}\subset\Lambda_{k}$) the $\Lambda_{k-1}%
$--subalgebra generated by $1$ and $\mathcal{C}\Lambda_{k}^{1}$
(resp.~$\mathcal{H}\Lambda_{k}^{1}$). Clearly, $d_{k}^{\,\mathsf{v}%
}(\mathcal{C}^{\bullet}\Lambda_{k})\subset\mathcal{C}^{\bullet}\Lambda_{k}$
and $d_{k}^{\,\mathsf{h}}(\mathcal{H}\Lambda_{k})\subset\mathcal{H}\Lambda
_{k}$. Moreover, $\Lambda_{k}^{1}=\mathcal{C}\Lambda_{k}^{1}\oplus
\mathcal{H}\Lambda_{k}^{1}$ and $\Lambda_{k}\simeq\mathcal{C}^{\bullet}%
\Lambda_{k}\otimes_{\Lambda_{k-1}}\mathcal{H}\Lambda_{k}$.

Put $\mathrm{D}^{\,\mathsf{v}}(\Lambda_{k-1})\overset{\mathrm{def}}{=}%
\{X\in\mathrm{D}(\Lambda_{k-1},\Lambda_{k-1})\;|\;X\circ\pi_{\infty}^{\ast
}=0\}$. Elements of $\mathrm{D}^{\,\mathsf{v}}(\Lambda_{k-1})$ are
called
vertical derivations of $\Lambda_{k-1}$. $\mathrm{D}(\Lambda_{k-1}%
,\Lambda_{k-1})$ splits as $\mathrm{D}(\Lambda_{k-1},\Lambda_{k-1}%
)\simeq\mathcal{C}\mathrm{D}(\Lambda_{k-1})\oplus\mathrm{D}^{\,\mathsf{v}%
}(\Lambda_{k-1})$ (see \cite{vv07}). Moreover, $\mathcal{C}\mathrm{D}%
(\Lambda_{k-1})\simeq\mathrm{Hom}_{\Lambda_{k-1}}(\mathcal{H}\Lambda_{k}%
^{1},\Lambda_{k-1})$ and $\mathrm{D}^{\,\mathsf{v}}(\Lambda_{k-1}%
)\simeq\mathrm{Hom}_{\Lambda_{k-1}}(\mathcal{C}\Lambda_{k}^{1},\Lambda_{k-1}%
)$. Denote by
$V_{j}^{\sigma,K}\in\mathrm{D}^{\,\mathsf{v}}(\Lambda_{k-1})$ and
$D_{\mu}^{K}\in\mathcal{C}\mathrm{D}(\Lambda_{k-1})$ with
$K\subset \{1,\ldots,k-1\}$, $j=1,\ldots,m$, $\mu=1,\ldots,n$ and
$\sigma$ being a lenght $n$ multi--index, elements of the dual local bases of
(\ref{basis1}) and (\ref{basis2}), respectively.

The following two sub--algebras of $\Lambda_{k-1}$ will be of a
special interest in the sequel.

\begin{itemize}
\item $\mathcal{C}_{\star}\Lambda_{k-1}$: the $\mathcal{F}$--subalgebra
generated by elements of the form $d_{K}^{\,\mathsf{v}}f$,
$f\in\mathcal{F}$,
$K\subset\{1,\ldots,k-1\}$ (note that $\mathcal{C}_{\star}\Lambda_{0}%
\equiv\mathcal{C}_{\star}\Lambda=\mathcal{C}^{\bullet}\Lambda$).

\item $\mathcal{C}_{\circ}\Lambda_{k-1}$: the $C^{\infty}(E)$--subalgebra
generated by elements of the form $d_{K}^{\,\mathsf{v}}g$, $g\in
C^{\infty }(E)$, $K\subset\{1,\ldots,k-1\}$.
\end{itemize}

Clearly, $\mathcal{C}_{\circ}\Lambda_{k-1}\subset\mathcal{C}_{\star}%
\Lambda_{k-1}$. Moreover,
$V_{j}^{\sigma,K}(\mathcal{C}_{\circ}\Lambda
_{k-1})\subset\mathcal{C}_{\circ}\Lambda_{k-1}$ with
$K\subset\{1,\ldots ,k-1\},\;j=1,\ldots,m$ and $\sigma$ being a
lenght $n$ multi--index. Denote by $\Lambda _{k-1}\varkappa$ the
$\mathcal{C}_{\star}\Lambda_{k-1}$--module of
$\mathcal{C}_{\star}\Lambda_{k-1}$--valued derivations of
$\mathcal{C}_{\circ }\Lambda_{k-1}$. $\Lambda_{k-1}\varkappa$ is
locally generated by derivations $V_{j}^{K}\equiv
V_{j}^{(0,\ldots,0),K}$, $K\subset\{1,\ldots
,k-1\},\;j=1,\ldots,m$.

$k$--times IDFs--valued symmetries (see \cite{vv07}) of $J^{\infty}(\pi)$ are
easily described in terms of vertical derivations of $\Lambda_{k-1}$. Indeed,
let $\chi=[X]\in\Lambda_{k-1}\mathrm{Sym}$, $X\in\mathrm{D}_{\mathcal{C}%
}(\Lambda_{k-1})$. There is a unique vertical representative of
$\chi$ denoted by $\rE_{\chi}$. $\rE_{\chi}$ is called
an \emph{evolutionary
derivation} of $\Lambda_{k-1}$. $\rE_{\chi}(\mathcal{C}_{\star}\Lambda_{k-1})\subset\mathcal{C}%
_{\star}\Lambda_{k-1}$, since $\rE_{\chi}\in\mathrm{D}%
_{\mathcal{C}}(\Lambda_{k-1})\cap\mathrm{D}^{\,\mathsf{v}}(\Lambda_{k-1})$.

\begin{proposition}
The correspondence $\Lambda_{k-1}\mathrm{Sym}\ni\chi\longmapsto\rE_{\chi
}|_{\mathcal{C}_{\circ}\Lambda_{k-1}}\in\Lambda_{k-1}\varkappa$ is an
isomorphism of vector spaces, whose inverse looks locally as $\Lambda
_{k-1}\varkappa\ni\chi_{K}^{j}V_{j}^{K}\longmapsto(D_{\sigma}\chi_{K}%
^{j})V_{j}^{\sigma,K}\in\Lambda_{k-1}\mathrm{Sym}$ with $\chi_{K}^{j}%
\in\mathcal{C}_{\star}\Lambda_{k-1}$,
$D_{\sigma}=(D_{1}^{(0,\ldots,0) })^{\sigma_1}\circ\cdots\circ
(D_{n}^{(0,\ldots,0)})^{\sigma_n}$, $K\subset\{1,\ldots,k-1\}$,
$j=1,\ldots,m$ and $\sigma=(\sigma_{1},\ldots,\sigma_{n})$ being a lenght $n$
multi--index.
\end{proposition}

As a consequence $\Lambda_{k-1}\mathrm{Sym}$ inherits the
structure of a $\mathcal{C}_{\star}\Lambda_{k-1}$--module. On the
other hand $\Lambda _{k-1}\varkappa$ inherits the structure of a
graded Lie--algebra denoted by
$(\Lambda_{k-1}\varkappa,\{\cdot,\cdot\})$. In the following we
identify $\Lambda_{k-1}\mathrm{Sym}$ and $\Lambda_{k-1}\varkappa$.
Finally, note that $d_{m}^{\,\mathsf{v}}$ is an evolutionary
derivation of $\Lambda_{k-1}$ for any $m<k$. Denote by $U_{m}$ the
corresponding element in $\Lambda _{k-1}\varkappa$, i.e.,
$U_{m}=d_{m}^{\,\mathsf{v}}|_{\mathcal{C}_{\circ }\Lambda_{k-1}}$.

\section{Adjoint graded $\mathcal{C}$--differential operators and horizontal
modules}

In what follows the algebra $\mathcal{C}_{\star}\Lambda_{k-1}$ is
considered to be the ground algebra. In particular, all
differential operators are operators over this algebra. Obviously,
$D_{\sigma}(\mathcal{C}_{\star}\Lambda
_{k-1})\subset\mathcal{C}_{\star}\Lambda_{k-1}$ for any
lenght $n$ multi--index $\sigma$. Let $P,Q$ be locally free graded
$\mathcal{C}_{\star}\Lambda_{k-1}$--modules of finite rank.

\begin{definition}
A linear differential operator $\square:P\longrightarrow Q$ is
called $\mathcal{C}$--differential if for any local basis
$\{e_{1},\ldots,e_{r}\}$ of $P$, $\square$ is locally of the form
$\square(p)=(-1)^{|p|\cdot|\alpha
|}\square_{\alpha}^{\sigma}D_{\sigma}p^{\alpha}$, where
$\square_{\alpha }^{\sigma}\in Q$ and $p=p^{\alpha}e_{\alpha}\in
P$, $|\alpha|{}\equiv |e_{\alpha}|$,
$p^{\alpha}\in\mathcal{C}_{\star}\Lambda_{k-1}$, $\alpha
=1,\ldots,r$.
\end{definition}

The totality of all $\mathcal{C}$--differential operators $\square
:P\longrightarrow Q$ has a natural
$\mathcal{C}_{\star}\Lambda_{k-1}$--module structure. This module
is denoted by $\mathcal{C}\mathrm{Diff}(P,Q)$. Similarly, denote
by $\mathcal{C}\mathrm{Diff}_{(p)}^{\mathrm{alt}}(P,Q)$ the
$\mathcal{C}_{\star
}\Lambda_{k-1}$--module of $Q$--valued, skew--symmetric, multi--$\mathcal{C}%
$--differential operators over $P$ with $p\geq0$ entries ($\mathcal{C}%
\mathrm{Diff}_{(0)}^{\mathrm{alt}}(P,Q)\equiv Q$). $\mathcal{H}\Lambda_{k}$ is
a $\mathcal{C}_{\star}\Lambda_{k-1}$--module and $d_{k}^{\,\mathsf{h}%
}:\mathcal{H}\Lambda_{k}\longrightarrow\mathcal{H}\Lambda_{k}$ is a
$\mathcal{C}$--differential operator.

The theory of adjoint $\mathcal{C}$--differential operators in the
category of $\mathcal{C}_{\star}\Lambda_{k-1}$--modules may be
developed almost literally as in the standard non--graded case
(see \cite{b99}). Here we limit ourselves to those elements of this
theory that are necessary for the proposed below description of
the $\Lambda_{k-1}\mathcal{C}$--spectral sequence terms .

Let $P$ be a graded $\mathcal{C}_{\star}\Lambda_{k-1}$--module. Consider the
map $w_{k}^{P}:\mathcal{C}\mathrm{Diff}(P,\mathcal{H}\Lambda_{k}%
)\longrightarrow\mathcal{C}\mathrm{Diff}(P,\mathcal{H}\Lambda_{k})$,
given by $w_{k}^{P}(\square)=d_{k}^{\,\mathsf{h}}\circ\square$,
$\square\in \mathcal{C}\mathrm{Diff}(P,\mathcal{H}\Lambda_{k})$.
$w_{k}^{P}$ is a differential, i.e., $w_{k}^{P}\circ w_{k}^{P}=0$.
The cohomology space of $w_{k}^{P}$ carries a natural
$\mathcal{C}_{\star}\Lambda_{k-1}$--module structure defined by
$\omega\cdot\lbrack
\square]\overset{\mathrm{def}}{=}(-1)^{|\omega|\cdot|\square|}[\square
\circ\omega]$, $\omega\in\mathcal{C}_{\star}\Lambda_{k-1}$,
$\square
\in\mathcal{C}\mathrm{Diff}(P,\mathcal{H}\Lambda_{k})$, $d_{k}^{\,\mathsf{h}%
}\circ\square=0$. The corresponding
$\mathcal{C}_{\star}\Lambda_{k-1}$--module is denoted by
$\widehat{P}$ and called the \emph{adjoint to $P$ module}. Put
$\Lambda_{k-1}\mathcal{B}\equiv\widehat{\mathcal{C}_{\star}\Lambda_{k-1}}$.

\begin{proposition}
\label{Prop}$\Lambda_{k-1}\mathcal{B}\simeq\mathcal{C}_{\star}\Lambda
_{k-1}\otimes_{\mathcal{F}}\mathcal{H}\Lambda_{1}^{n}$. Moreover,
if $P$ is a locally free
$\mathcal{C}_{\star}\Lambda_{k-1}$--module of finite rank, then
$\widehat{P}\simeq\mathrm{Hom}_{\mathcal{C}_{\star}\Lambda_{k-1}}(P,\Lambda_{k-1}\mathcal{B})$.
\end{proposition}

In particular, proposition \ref{Prop} tells that
$\Lambda_{k-1}\mathcal{B}$ is a locally free
$\mathcal{C}_{\star}\Lambda_{k-1}$--module of rank one. Thus, if
$P$ is a locally free $\mathcal{C}_{\star}\Lambda_{k-1}$--module
of finite rank, then $\widehat{\widehat{P}}\simeq P$. Further on
our exposition is based on identifications of proposition
\ref{Prop}. Note that $\Lambda _{k-1}\mathcal{B}$ is embedded in
$\ker d_{k}^{\,\mathsf{h}}\subset \mathcal{H}\Lambda_{k}$ in a
natural way by means of the correspondence $\Lambda
_{k-1}\mathcal{B}\simeq\mathcal{C}_{\star}\Lambda_{k-1}\otimes_{\mathcal{F}%
}\mathcal{H}\Lambda_{1}^{n}\ni\omega\otimes\sigma\mapsto\omega
\cdot\kappa_{1k}(\sigma)\in\mathcal{H}\Lambda_{k}$, $\omega\in\mathcal{C}%
_{\star}\Lambda_{k-1}$, $\sigma\in\mathcal{H}\Lambda_{1}^{n}$.

If $P,Q$ are $\mathcal{C}_{\star}\Lambda_{k-1}$--modules and
$\square \in\mathcal{C}\mathrm{Diff}(P,Q)$, then a differential
operator
$\widehat{\square}\in\mathcal{C}\mathrm{Diff}(\widehat{Q},\widehat
{P})$   is well--defined by
$\widehat{\square}[\Delta]=(-1)^{|\Delta|\cdot|\square|}[\Delta
\circ\square]$,
$\Delta\in\mathcal{C}\mathrm{Diff}(Q,\mathcal{H}\Lambda_{k})$,
$d_{k}^{\,\mathsf{h}}\circ\Delta=0$. $\widehat{\square}$ is called
the \emph{adjoint to $\square$ operator}.

Let $P$ be as above, $\xi_{1},\ldots,\xi_{p-2}\in P$ and $\square
\in\mathcal{C}\mathrm{Diff}_{(p-1)}^{\mathrm{alt}}(P,\widehat{P})$. Define a
$\mathcal{C}$--differential operator $\square_{\xi_{1},\ldots,\xi_{p-2}%
}:P\longrightarrow\widehat{P}$ by putting $\square_{\xi_{1},\ldots,\xi_{p-2}%
}(\xi)\overset{\mathrm{def}}{=}\square(\xi_{1},\ldots,\xi_{p-2},\xi
)\in\widehat{P}$. Also put
\[
L_{p}^{(k)}(P)\overset{\mathrm{def}}{=}\{\square\in\mathcal{C}\mathrm{Diff}%
_{(p-1)}^{\mathrm{alt}}(P,\widehat{P})\;|\;\widehat{\square_{\xi_{1}%
,\ldots,\xi_{p-2}}}=-\square_{\xi_{1},\ldots,\xi_{p-2}},\;\forall\xi
_{1},\ldots,\xi_{p-2}\in P\}.
\]

A $\mathcal{C}_{\star}\Lambda_{k-1}$--module $P$ is called
\emph{horizontal} if it is
of the form $P\simeq\mathcal{C}_{\star}\Lambda_{k-1}\otimes_{C^{\infty}%
(M)}P_{0}$ for a (possibly graded) $C^{\infty}(M)$--module
$P_{0}$. Evolutionary derivations of $\Lambda_{k-1}$ act naturally
on horizontal modules. Namely, let
$\chi\in\Lambda_{k-1}\mathrm{Sym}$ and $\rE_{\chi}$ be the
corresponding evolutionary derivation of $\Lambda_{k-1}$. For $\xi
=\omega\otimes\xi_{0}\in P$,
$\omega\in\mathcal{C}_{\star}\Lambda_{k-1}$ and $\xi_{0}\in P_{0}$
we put
\begin{equation}
\rE_{\chi}\xi\overset{\mathrm{def}}{=}\rE_{\chi}\omega\otimes
\xi_{0}\in P.\label{DefHor}%
\end{equation}
Since $\rE_{\chi}$ is a vertical derivation, definition
(\ref{DefHor}) extends unambiguously to the whole $P$. Now, fix an
element $\xi\in P$ and define an operator
$\ell_{\xi}^{\{k\}}:\Lambda_{k-1}\varkappa
\longrightarrow P$ by putting $\ell_{\xi}^{\{k\}}(\chi)\overset{\mathrm{def}%
}{=}(-1)^{|\xi|\cdot|\chi|}\rE_{\chi}\xi$, $\chi\in\Lambda
_{k-1}\varkappa$. $\ell_{\xi}^{\{k\}}$ is a
$\mathcal{C}$--differential operator called the \emph{universal
}$\Lambda_{k-1}$\emph{--linearization} of $\xi$.

\begin{example}
$\mathcal{H}\Lambda_{k}$ is a horizontal module. Indeed, $\mathcal{H}%
\Lambda_{k}\simeq\mathcal{C}_{\star}\Lambda_{k-1}\otimes_{C^{\infty}%
(M)}\Lambda_{k}(M)$.
\end{example}

\begin{example}
$\Lambda_{k-1}\varkappa$ and $\widehat{\Lambda_{k-1}\varkappa}$ are horizontal
modules. Indeed, let $W_{k-1}\subset\Lambda_{k-1}\varkappa$ be the $C^{\infty
}(M)$--submodule of derivations that are locally of the form $\chi_{K}%
^{j}V_{j}^{K}$ with $\chi_{K}^{j}\in C^{\infty}(M)$,
$K\subset\{1,\ldots,k-1\}$,
$j=1,\ldots,m$. $W_{k-1}$ is well defined. Then $\Lambda_{k-1}%
\varkappa\simeq\mathcal{C}_{\star}\Lambda_{k-1}\otimes_{C^{\infty}(M)}W_{k-1}$
and $\widehat{\Lambda_{k-1}\varkappa}\simeq\mathcal{C}_{\star}\Lambda
_{k-1}\otimes_{C^{\infty}(M)}\mathrm{Hom}_{C^{\infty}(M)}(W_{k-1},\Lambda
^{n}(M))$.
\end{example}

\section{Secondary IDF\lowercase{s} on $J^{\infty}(\pi)$}

The complexes $(\Lambda_{k-1}\mathcal{C}E_{0}^{0,\bullet}=\Lambda
_{k}/\mathcal{C}\Lambda_{k},d_{k,0}^{0,\bullet})$ and
$(\mathcal{H}\Lambda _{k},d_{k}^{\,\mathsf{h}})$ are
isomorphic in a natural way and further on
we shall identify $\Lambda_{k-1}\mathcal{C}%
E_{0}^{0,\bullet}$ with $\mathcal{H}\Lambda_{k}$.

\begin{remark}
For any $p>0$, there exists a natural isomorphism of complexes
\[
\xymatrix@R=30pt{ \Lambda_{k-1}\mathcal{C}E_0^{p,\bullet}
\ar[r]^-{d_{0,k}^{p,\bullet}}
\ar[d]_{\eta_{k-1}}^{\begin{sideways}$\widetilde{\quad\quad}$\end{sideways}}
& \Lambda_{k-1}\mathcal{C}E_0^{p,\bullet}
\ar[d]_{\eta_{k-1}}^{\begin{sideways}
$\widetilde{\quad\quad}$\end{sideways}} \\
\mathcal{C}\mathrm{Diff}^{\mathrm{alt}}_{(p)}(\Lambda_{k-1}\varkappa,
\mathcal{H}\Lambda_k) \ar[r]^{w^p_k} &
\mathcal{C}\mathrm{Diff}^{\mathrm{alt}}_{(p)}
(\Lambda_{k-1}\varkappa,\mathcal{H}\Lambda_k) } ,
\]
where $w_{k}^{p}(\square)\overset{\mathrm{def}}{=}(-1)^{p}d^{\,\mathsf{h}%
}\circ\square$, $\square\in\mathcal{C}\mathrm{Diff}_{(p)}^{\mathrm{alt}%
}(\Lambda_{k-1}\varkappa,\mathcal{H}\Lambda_{k})$. The isomorphism
$\eta_{k-1}$ is defined by
\[
\eta_{k-1}([\omega]_{\mathcal{C}^{p+1}\Lambda_{k}})(\chi_{1},\ldots,\chi
_{p})=(-1)^{|\omega|\cdot(|\chi_{1}|+\cdots+|\chi_{p}|)+p(p-1)/2}%
[(i_{\rE_{\chi_{1}}}^{\{k\}}\circ\cdots\circ i_{\rE_{\chi_{p}%
}}^{\{k\}})(\omega)]_{\mathcal{C}\Lambda_{k}}\in\mathcal{H}\Lambda_{k}.
\]

\end{remark}

The following two theorems generalize the results of
\cite{v78,v84} concerning the first term of the
$\mathcal{C}$--spectral sequence to that of the
$\Lambda_{k-1}\mathcal{C}$--spectral one. In particular, they give
an explicit description of secondary IDFs on $J^{\infty}(\pi)$ and
basic operations over them.

\begin{theorem}
[One Line Theorem]\quad

\begin{itemize}
\item $\Lambda_{k}\mathcal{C}E_{1}^{p,q}=0$ if either $q>n$, or $p>0$ and $q<n$;

\item $\Lambda_{k-1}\mathcal{C}E_{1}^{0,\bullet}\simeq H(\mathcal{H}%
\Lambda_{k},d_{k}^{\,\mathsf{h}})\simeq\Lambda_{k-2}\mathcal{C}E_{1}%
^{\bullet,\bullet}$;

\item $\Lambda_{k-1}\mathcal{C}E_{1}^{p,n}\simeq L_{p}^{(k)}(\Lambda
_{k-1}\varkappa)$, if $p>0$.
\end{itemize}
\end{theorem}

\begin{corollary}
$\Lambda_{k-1}\mathcal{C}E$ stabilizes at the second term and
$\Lambda _{k-1}\mathcal{C}E_{2}^{0,q}\simeq
H^q(\Lambda_{k},d_{k})\simeq H^q(E), \;q\leq n$, and
$\Lambda_{k-1}\mathcal{C}E_{2}^{p,n}=H^{p+n}(E), \;p\geq 0$.
\end{corollary}

There is a distinguished element $\mathfrak{u}$ in $\widehat{\mathcal{H}%
\Lambda_{k}}\simeq\mathrm{Hom}_{\mathcal{C}_{\star}\Lambda_{k-1}}%
(\mathcal{H}\Lambda_{k},\Lambda_{k-1}\mathcal{B})$ whose local
expression is
\begin{align*}
\mathfrak{u}(\Omega_{\mu_{1}\cdots\mu_{r}}^{J_{1}\cdots J_{r}}d_{J_{1}}%
x^{\mu_{1}}\cdots d_{J_{r}}x^{\mu_{r}})\overset{\mathrm{def}}{=} &
\Omega_{\mu_{1}\cdots\mu_{r}}^{J_{1}\cdots J_{r}}\otimes((\kappa_{1k}\circ
\nu)(d_{J_{1}}x^{\mu_{1}}\cdots d_{J_{r}}x^{\mu_{r}}))\\
{}\in{} &  \mathcal{C}_{\star}\Lambda_{k-1}\otimes_{\mathcal{F}}%
\mathcal{H}\Lambda_{1}^{n}\simeq\Lambda_{k-1}\mathcal{B},
\end{align*}
where $\Omega_{\mu_{1}\cdots\mu_{r}}^{J_{1}\cdots J_{r}}\in\mathcal{C}_{\star
}\Lambda_{k-1}$ and $\nu:\mathcal{H}\Lambda_{k}\longrightarrow\mathcal{H}%
\Lambda_{k}^{n}$ is the projection onto the homogeneous component
of multi--degree $(0,\ldots,0,n)\in%
%TCIMACRO{\U{2124} }%
%BeginExpansion
\mathbb{Z}
%EndExpansion
^{k}$.

\begin{theorem}
\quad

\begin{itemize}
\item The differential $d_{k,1}^{0,n}:H^{n}(\mathcal{H}\Lambda_{k}%
,d_{k}^{\,\mathsf{h}})\longrightarrow\widehat{\Lambda_{k-1}\varkappa}$ is
given by
\[
d_{k,1}^{0,n}[\Omega]_{\mathrm{im\,}d_{k}^{\,\mathsf{h}}}=\widehat{\ell
}_{\Omega}^{\,\{k\}}(\mathfrak{u}),\quad\Omega\in\mathcal{H}\Lambda_{k}^{n};
\]
and its local description is
\begin{align*}
d_{k,1}^{0,n}([\Omega]_{\mathrm{im\,}d_{k}^{\,\mathsf{h}}})(\chi) &
=(-1)^{\sigma+|\chi|\cdot|\Omega|}\chi_{L}^{j}(D_{\sigma}\circ V_{j}^{\sigma
L})(A)\otimes dx^{1}\cdots dx^{n}\\
{} &  \in{}\mathcal{C}_{\star}\Lambda_{k-1}\otimes_{\mathcal{F}}%
\mathcal{H}\Lambda_{1}^{n}\simeq\Lambda_{k-1}\mathcal{B},
\end{align*}
assuming that $\Omega=A\,d_{k}x^{1}\cdots
d_{k}x^{n}\in\mathcal{H}\Lambda_{k}^{n}$,
$A\in\mathcal{C}_{\star}\Lambda_{k-1}$ and $\chi=\chi_{L}^{j}V_{j}^{L}%
\in\Lambda_{k-1}\varkappa$;

\item the differential $d_{k,1}^{p,n}:L_{p}^{(k)}(\Lambda_{k-1}\varkappa
)\longrightarrow L_{p+1}^{(k)}(\Lambda_{k-1}\varkappa)$ acting on
secondary iterated $p$--forms, $p>0$, is given by

\begin{align*}
& d_{k,1}^{p,n}(\Theta)(\chi_{1},\ldots,\chi_{p})=\sum_{i=1}^{p}%
(-1)^{a(i)+i+1}(\rE_{\chi_{i}})(\Theta(\chi_{1},\ldots,\widehat{\chi_{i}%
},\ldots,\chi_{p}))\\
& +\sum_{i<j}(-1)^{c(i,j)+i+j}\Theta(\{\chi_{i},\chi_{j}\},\chi_{1}%
,\ldots,\widehat{\chi_{i}},\ldots,\widehat{\chi_{j}},\ldots,\chi_{p})\\
& +\tfrac{1}{p}\sum_{i=1}^{p}(-1)^{i+1}[(-1)^{a(i)}(p-1)\widehat{\ell}%
_{\chi_{i}}^{\,\{k+1\}}(\Theta(\chi_{1},\ldots,\widehat{\chi_{i}},\ldots
,\chi_{p}))\\
& -(-1)^{b(i)}\widehat{\ell}_{\Theta(\chi_{1},\ldots,\widehat{\chi_{i}}%
,\ldots,\chi_{p})}^{\,\{k+1\}}(\chi_{i})],\quad
\end{align*}
$\Theta\in L_{p}^{(k)}(\Lambda_{k-1}\varkappa)$, where $a(i)=|\chi_{i}|\cdot(|\Theta|+|\chi_{1}|+\cdots+|\chi_{i-1}|)$,
$b(i)=|\chi_{i}|\cdot(|\chi_{i+1}|+\cdots+|\chi_{p}|)$,
$c(i,j)=a(i)+a(j)-|\Theta|\cdot(|\chi_{i}|+|\chi_{j}|)-|\chi_{i}|\cdot
|\chi_{j}|$;

\item the action $\mathcal{L}_{\chi}^{\{k\}}:L_{p}^{(k)}(\Lambda
_{k-1}\varkappa)\longrightarrow L_{p}^{(k)}(\Lambda_{k-1}\varkappa)$ of
$\chi\in\Lambda_{k-1}\mathrm{Sym}$ on secondary iterated $p$--forms, $p>0$, is
given by%

\begin{align*}
(\mathcal{L}_{\chi}^{\{k\}}\Theta)(\chi_{1},\ldots,\chi_{p-1})(\chi_{p})  &
=(\rE_{\chi}\circ\Theta)(\chi_{1},\ldots,\chi_{p})\\
& -\sum_{i=1}^{p}(-1)^{|\chi|\cdot(|\chi_{1}|+\cdots+|\chi_{i-1}|+|\Theta
|)}\Theta(\chi_{1},\ldots,\{\chi,\chi_{i}\},\ldots,\chi_{p}),
\end{align*}
$\Theta\in L_{p}^{(k)}(\Lambda_{k-1}\varkappa)$.

\item the insertion $i_{\chi}^{\{k\}}:\widehat{\Lambda_{k-1}\varkappa
}\longrightarrow H(\mathcal{H}\Lambda_{k},d_{k}^{\,\mathsf{h}})$ of $\chi
\in\Lambda_{k-1}\mathrm{Sym}$ in secondary iterated $1$--forms is given by
\[
i_{\chi}^{\{k\}}\Theta=(-1)^{|\chi|\cdot|\Theta|}[\Theta(\chi)]_{\mathrm{im\,}%
d_{k}^{\,\mathsf{h}}},\quad\Theta\in\widehat{\Lambda_{k-1}\varkappa}%
,\;\Theta(\chi)\in\Lambda_{k-1}\mathcal{B}\hookrightarrow\ker d_{k}%
^{\,\mathsf{h}}.
\]

\item the insertion $i_{\chi}^{\{k\}}:L_{p}^{(k)}(\Lambda_{k-1}\varkappa
)\longrightarrow L_{p-1}^{(k)}(\Lambda_{k-1}\varkappa)$ of $\chi\in
\Lambda_{k-1}\mathrm{Sym}$ in secondary iterated $p$--forms, $p>1$, is given
by%
\[
(i_{\chi}^{\{k\}}\Theta)(\chi_{1},\ldots,\chi_{p-2})=(-1)^{|\chi|\cdot
|\Theta|}\Theta(\chi,\chi_{1},\ldots,\chi_{p-2}).
\]

\end{itemize}
\end{theorem}

\section{Secondary covariant tensors}

Recall that $\tau\in\Lambda_{k-1}\mathcal{C}E_{1}^{(1,\ldots,1),\bullet}$ is a
secondary covariant $k$--tensor iff $\mathcal{L}_{I_{m}^{K^{\prime}}}%
^{\{k\}}\tau=i_{I_{m}^{K}}^{\{k\}}\tau=0$ for any $m<k$ and $K=\{k_{1}%
,\ldots,k_{s}\},K^{\prime}=\{k_{1}^{\prime},\ldots,k_{s^{\prime}}^{\prime
}\}\subset\{1,\ldots,k-1\}$, $s\geq2$, $s^{\prime}\geq1$ (see \cite{vv07}). In
what follows we characterize such $\tau$'s. First of all, there is a well
defined $\mathcal{F}$--linear surjective map $\mathfrak{p}_{k-1}:\Lambda
_{k-1}\varkappa\longrightarrow\Lambda_{0}\varkappa\equiv\varkappa$ locally
given by
\[
\mathfrak{p}_{k-1}(\chi_{K}^{j}V_{j}^{K})\overset{\mathrm{def}}{=}%
\mathsf{p}_{0}(\chi^{j})\tfrac{\partial}{\partial u^{j}},\quad\chi_{K}^{j}%
\in\mathcal{C}_{\star}\Lambda_{k-1},\;K\subset\{1,\ldots,k-1\},\;j=1,\ldots,m,
\]
$\mathsf{p}_{0}:\mathcal{C}_{\star}\Lambda_{k-1}\longrightarrow\mathcal{F}$
being the projection onto the  homogeneous component of multi--degree $(0,\ldots,0)\in%
%TCIMACRO{\U{2124} }%
%BeginExpansion
\mathbb{Z}
%EndExpansion
^{k-1}$. Then there is an injective morphism of
$\mathcal{F}$--modules
$\mathfrak{i}_{k}:(\mathcal{C}\Lambda^{1})^{\otimes
k-1}\otimes_{\mathcal{F}}\widehat{\varkappa}\longrightarrow\Lambda
_{k}\mathcal{C}E_{1}^{(1,\ldots,1),\bullet}\subset\widehat{\Lambda
_{k-1}\varkappa}$ given by
\begin{align*}
\mathfrak{i}_{k}(\omega_{1}\otimes\cdots\otimes\omega_{k-1}\otimes\psi
)(\chi)\overset{\mathrm{def}}{=} &  (\omega_{1}\cdot\kappa_{12}(\omega
_{2})\cdot\cdots\cdot\kappa_{1k-1}(\omega_{k-1}))\otimes((\psi\circ
\mathfrak{p}_{k-1})(\chi))\\
{}\in{} &  \mathcal{C}_{\star}\Lambda_{k-1}\otimes_{\mathcal{F}}%
\mathcal{H}\Lambda_{1}^{n}\simeq\Lambda_{k-1}\mathcal{B},
\end{align*}
$\omega_{1},\ldots,\omega_{k-1}\in\mathcal{C}\Lambda^{1}$, $\psi\in
\widehat{\varkappa}$, $\chi\in\Lambda_{k-1}\varkappa$.

Now we are able to characterize secondary IDFs on
$J^{\infty}(\pi)$ in the same manner as it was done for ordinary
IDFs in \cite{vv06}
\begin{proposition}
$\tau\in\Lambda_{k-1}\mathcal{C}E_{1}^{(1,\ldots,1),\bullet}$ is a secondary
covariant $k$--tensor iff $\tau\in\mathrm{im\,}\mathfrak{i}_{k}$.
\end{proposition}

\section{Conclusions}

The above given description of secondary IDFs on $J^{\infty}(\pi)$
was not conceptual in all its aspects just because the necessary
for this purpose constructions can not be put in frames of a short
note. This will be done in the subsequent detailed exposition
together with basic elements of secondary tensor calculus and some
applications to mechanics and field theory.

\end{document}